\documentclass[12pt]{article}%
\usepackage{amsfonts}
\usepackage{amssymb}
\usepackage{amsmath}
\usepackage{amscd}
\usepackage{latexsym}
\usepackage[onehalfspacing]{setspace}
\usepackage{theorem}
\usepackage{natbib}
\usepackage[colorlinks=true,urlcolor=black,linkcolor=black,pdfpagemode="None",pdfstartview=FitH]%
{hyperref}
\usepackage{graphicx}
\usepackage{caption}
\usepackage{subfig}
\usepackage{lscape}%
\providecommand{\U}[1]{\protect\rule{.1in}{.1in}}
\newtheorem{theorem}{Theorem}

\newtheorem{assumption}{Assumption}
\newcommand{\I}{{\rm 1\hspace*{-0.4ex}\rule{0.1ex}{1.52ex}\hspace*{0.2ex}}}
\hypersetup{pdfborder = {0 0 0},colorlinks=true,linkcolor=blue,citecolor=blue}
\renewcommand{\cite}{\citet*}
\begin{document}

\title{\vspace{-0.6in}Inference in Linear Regression Models with Many Covariates and
Heteroskedasticity\thanks{We thank Xinwei Ma, Ulrich M\"{u}ller and Andres
Santos for very thoughtful discussions regarding this project. We also thank
Silvia Gon\c{c}alvez, Pat Kline and James MacKinnon. In addition, an Associate
Editor and three reviewers offered excellent recommendations that improved
this paper. The first author gratefully acknowledges financial support from
the National Science Foundation (SES 1459931). The second author gratefully
acknowledges financial support from the National Science Foundation (SES
1459967) and the research support of CREATES (funded by the Danish National
Research Foundation under grant no. DNRF78).}\bigskip}
\author{Matias D. Cattaneo\thanks{Department of Economics and Department of
Statistics, University of Michigan.}
\and Michael Jansson\thanks{Department of Economics, UC Berkeley and \emph{CREATES}%
.}
\and Whitney K. Newey\thanks{Department of Economics, MIT.}}
\maketitle

\begin{abstract}
The linear regression model is widely used in empirical work in Economics,
Statistics, and many other disciplines. Researchers often include many
covariates in their linear model specification in an attempt to control for
confounders. We give inference methods that allow for many covariates and
heteroskedasticity. Our results are obtained using high-dimensional
approximations, where the number of included covariates are allowed to grow as
fast as the sample size. We find that all of the usual versions of
Eicker-White heteroskedasticity consistent standard error estimators for
linear models are inconsistent under this asymptotics. We then propose a new
heteroskedasticity consistent standard error formula that is fully automatic
and robust to both (conditional)\ heteroskedasticity of unknown form and the
inclusion of possibly many covariates. We apply our findings to three
settings: parametric linear models with many covariates, linear panel models
with many fixed effects, and semiparametric semi-linear models with many
technical regressors. Simulation evidence consistent with our theoretical
results is also provided. The proposed methods are also illustrated with an
empirical application.\bigskip

\end{abstract}

\textit{Keywords:} high-dimensional models, linear regression, many
regressors, heteroskedasticity, standard errors.\bigskip%

\thispagestyle{empty}
\setcounter{page}{0}
\newpage\doublespacing

\section{Introduction\label{section:intro}}

A key goal in empirical work is to estimate the structural, causal, or
treatment effect of some variable on an outcome of interest, such as the
impact of a labor market policy on outcomes like earnings or employment. Since
many variables measuring policies or interventions are not exogenous,
researchers often employ observational methods to estimate their effects. One
important method is based on assuming that the variable of interest can be
taken as exogenous after controlling for a sufficiently large set of other
factors or covariates.\ A major problem that empirical researchers face when
employing selection-on-observables methods to estimate structural effects is
the availability of many potential covariates. This problem has become even
more pronounced in recent years because of the widespread availability of
large (or high-dimensional) new data sets.

Not only it is often the case that substantive discipline-specific theory (or
intuition) will suggest a large set of variables that might be important, but
also researchers usually prefer to include additional \textquotedblleft
technical\textquotedblright\ controls constructed using indicator variables,
interactions, and other non-linear transformations of those variables.
Therefore, many empirical studies include very many covariates in order to
control for as broad array of confounders as possible. For example, it is
common practice to include dummy variables for many potentially overlapping
groups based on age, cohort, geographic location, etc. Even when some controls
are dropped after valid covariate selection
(\cite{Belloni-Chernozhukov-Hansen_2014_ReStud}), many controls usually may
remain in the final model specification. For example,
\cite{Angrist-Hahn_2004_ReStat} discuss when to include many covariates in
treatment effect models.

We present valid inference methods that explicitly account for the presence of
possibly many controls in linear regression models under
(conditional)\ heteroskedasticity. We consider the setting where the object of
interest is $\mathbf{\beta}$ in a model of the form%
\begin{equation}
y_{i,n}=\mathbf{\beta}^{\prime}\mathbf{x}_{i,n}+\mathbf{\gamma}_{n}^{\prime
}\mathbf{w}_{i,n}+u_{i,n},\qquad i=1,\ldots,n,\label{eq:OLSmodel}%
\end{equation}
where $y_{i,n}$ is a scalar outcome variable, $\mathbf{x}_{i,n}$ is a
regressor of small\ (i.e., fixed) dimension $d$, $\mathbf{w}_{i,n}$ is a
vector of covariates of possibly \textquotedblleft large\textquotedblright%
\ dimension $K_{n}$, and $u_{i,n}$ is an unobserved error term. Two important
cases discussed in more detail below, are \textquotedblleft
flexible\textquotedblright\ parametric modeling of controls via basis
expansions such as higher-order powers and interactions (i.e., a series-based
formulation of the partially linear regression model), and models with many
dummy variables such as multi-way fixed effects and interactions thereof in
panel data models. In both cases conducting OLS-based inference on
$\mathbf{\beta}$ in (\ref{eq:OLSmodel}) is straightforward when the error
$u_{i,n}$ is homoskedastic and/or the dimension $K_{n}$ of the nuisance
covariates is modeled as a vanishing fraction of the sample size. The latter
modeling assumption, however, is inappropriate in applications with many dummy
variables and does not deliver a good distributional approximation when many
covariates are included.

Motivated by the above observations, this paper studies the consequences of
allowing the error $u_{i,n}$ in (\ref{eq:OLSmodel}) to be (conditionally)
heteroskedastic in a setting where the covariate $\mathbf{w}_{i,n}$ is
permitted to be high-dimensional in the sense that $K_{n}$ is allowed, but not
required, to\ be a non-vanishing fraction of the sample size. Our main purpose
is to investigate the possibility of constructing
heteroskedasticity-consistent variance estimators for the OLS estimator of
$\mathbf{\beta}$ in (\ref{eq:OLSmodel}) without (necessarily) assuming any
special structure on the part of the covariate $\mathbf{w}_{i,n}$. We present
two main results. First, we provide high-level sufficient conditions
guaranteeing a valid Gaussian distributional approximation to the finite
sample distribution of the OLS estimator of $\mathbf{\beta}$, allowing for the
dimension of the nuisance covariates to be \textquotedblleft
large\textquotedblright\ relative to the sample size (i.e., $K_{n}%
/n\not \rightarrow 0$). Second, we characterize the large sample properties of
a class of variance estimators, and use this characterization to obtain both
negative and positive results. The negative finding is that the Eicker-White
estimator is inconsistent in general, as are popular variants of this
estimator. The positive result gives conditions under which an alternative
heteroskedasticity-robust variance estimator (described in more detail below)
is consistent. The main condition needed for our constructive results is a
high-level assumption on the nuisance covariates requiring in particular that
their number be strictly less than half of the sample size. As a by-product,
we also find that among the popular HC$k$ class of standard errors estimators
for linear models, a variant of the HC$3$ estimator delivers standard errors
that are asymptotically upward biased in general. Thus, standard
OLS\ inference employing HC$3$ standard errors will be asymptotically valid,
albeit conservative, even in high-dimensional settings where the number of
covariate $\mathbf{w}_{i,n}$ is large relative to the sample size, i.e., when
$K_{n}/n\not \rightarrow 0$.

Our results contribute to the already sizeable literature on
heteroskedasticity-robust variance estimators for linear regression models, a
recent review of which is given by \cite{MacKinnon_2012_BookCh}. Important
papers whose results are related to ours include \cite{White_1980_ECMA},
\cite{MacKinnon-White_1985_JoE}, \cite{Wu_1986_AoS},
\cite{Chesher-Jewitt_1987_ECMA}, \cite{Shao-Wu_1987_AoS},
\cite{Chesher_1989_ECMA}, \cite{CribariNeto-Lima_2000_Biometrika},
\cite{Kauermann-Carroll_2001_JASA},
\cite{Bera-Suprayitno-Premaratne_2002_JSPI}, \cite{Stock-Watson_2008_ECMA},
\cite{CribariNeto-Lima_2011_JSPI}, \cite{Muller_2013_ECMA}, and
\cite{Abadie-Imbens-Zheng_2014_JASA}. In particular,
\cite{Bera-Suprayitno-Premaratne_2002_JSPI} analyze some finite sample
properties of a variance estimator similar to the one whose asymptotic
properties are studied herein. They use unbiasedness or minimum norm quadratic
unbiasedness to motivate a variance estimator that is similar in structure to
ours, but their results are obtained for fixed $K_{n}$ and $n$ and is silent
about the extent to which consistent variance estimation is even possible when
$K_{n}/n\not \rightarrow 0.$

This paper also adds to the literature on high-dimensional linear regression
where the number of regressors grow with the sample size; see, e.g.,
\cite{Huber_1973_AoS}, \cite{Koenker_1988_JAE}, \cite{Mammen_1993_AoS},
\cite{ElKaroui-Bean-Bickel-Lim-Yu_2013_PNAS},
\cite{Zheng-Jiang-Bai-He_2014_Biometrika}, \cite{Li-Muller_2017_wp}, and
references therein. In particular, \cite{Huber_1973_AoS} showed that fitted
regression values are not asymptotically normal when the number of regressors
grows as fast as sample size, while \cite{Mammen_1993_AoS} obtained asymptotic
normality for arbitrary contrasts of OLS estimators in linear regression
models where the dimension of the covariates is at most a vanishing fraction
of the sample size. More recently,
\cite{ElKaroui-Bean-Bickel-Lim-Yu_2013_PNAS} showed that, if a Gaussian
distributional assumption on regressors and homoskedasticity is assumed, then
certain estimated coefficients and contrasts in linear models are
asymptotically normal when the number of regressors grow as fast as sample
size, but do not discuss inference results (even under homoskedasticity). Our
result in Theorem \ref{thm:1} below shows that certain contrasts of
OLS\ estimators in high-dimensional linear models are asymptotically normal
under fairly general regularity conditions. Intuitively, we circumvent the
problems associated with the lack of asymptotic Gaussianity in general
high-dimensional linear models by focusing exclusively on a small subset of
regressors when the number of covariates gets large. We give inference results
by constructing heteroskedasticity consistent standard errors without imposing
any distributional assumption or other very specific restrictions on the
regressors. In particular, we do not require the coefficients $\mathbf{\gamma
}_{n}$ to be consistently estimated; in fact, they will not be in most of our
examples discussed below.

Our high-level conditions allow for $K_{n}\propto n$ and restrict the data
generating process in fairly general and intuitive ways. In particular, our
generic sufficient condition on the nuisance covariates $\mathbf{w}_{i,n}$
covers several special cases of interest for empirical work. For example, our
results encompass (and weakens in certain sense) those reported in
\cite{Stock-Watson_2008_ECMA}, who investigated the one-way fixed effects
panel data regression model and showed that the conventional Eicker-White
heteroskedasticity-robust variance estimator is inconsistent, being plagued by
a non-negligible bias problem attributable to the presence of many covariates
(i.e., the fixed effects). The very special structure of the covariates in the
one-way fixed effects estimator enables an explicit characterization of this
bias, and also leads to a direct plug-in consistent bias-corrected version of
the Eicker-White variance estimator. The generic variance estimator proposed
herein essentially reduces to this bias-corrected variance estimator in the
special case of the one-way fixed effects model, even though our results are
derived from a different perspective and generalize to other settings.

Furthermore, our general inference results can be used when many multi-way
fixed effects and similar discrete covariates are introduced in a linear
regression model, as it is usually the case in social interaction and network
settings. For example, in a very recent contribution, \cite{Verdier_2017_wp}
develops new results for two-way fixed effect design and projection matrices,
and use them to verify our high-level conditions in linear models with two-way
unobserved heterogeneity and sparsely matched data (which can also be
interpreted as a network setting). These results provide another interesting
and empirically relevant illustration of our generic theory.
\cite{Verdier_2017_wp} also develops inference results able to handle time
series dependence in his specific context, which are not covered by our
assumptions because we impose independence in the cross-sectional dimension of
the (possibly grouped)\ data.

The rest of this paper is organized as follows. Section \ref{section:overview}
presents the variance estimators we study and gives a heuristic description of
their main properties. Section \ref{section:framework} introduces our general
framework, discusses high-level assumptions and illustrates the applicability
of our methods using three leading examples. Section \ref{section:results}
gives the main results of the paper. Section \ref{section:simuls} reports the
results of a Monte Carlo experiment, while \ref{section:empapp} illustrates
our methods using an empirical application. Section \ref{section:conclusion}
concludes. Proofs and additional methodological and numerical results are
reported in the online supplemental appendix.

\section{Overview of Results\label{section:overview}}

For the purposes of discussing distribution theory and variance estimators
associated with the OLS estimator $\mathbf{\hat{\beta}}_{n}$ of $\mathbf{\beta
}$ in (\ref{eq:OLSmodel}), when possibly the $K_{n}$-dimensional nuisance
covariates $\mathbf{w}_{i,n}$ is of \textquotedblleft large\textquotedblright%
\ dimension and/or the parameters $\mathbf{\gamma}_{n}$ cannot be estimated
consistently, it is convenient to write the estimator in \textquotedblleft
partialled out\textquotedblright\ form as%
\[
\mathbf{\hat{\beta}}_{n}=(\sum_{i=1}^{n}\mathbf{\hat{v}}_{i,n}\mathbf{\hat{v}%
}_{i,n}^{\prime})^{-1}(\sum_{i=1}^{n}\mathbf{\hat{v}}_{i,n}y_{i,n}%
),\text{\qquad}\mathbf{\hat{v}}_{i,n}=\sum_{j=1}^{n}M_{ij,n}\mathbf{x}_{j,n},
\]
where $M_{ij,n}=%
\I
(i=j)-\mathbf{w}_{i,n}^{\prime}(\sum_{k=1}^{n}\mathbf{w}_{k,n}\mathbf{w}%
_{k,n}^{\prime})^{-1}\mathbf{w}_{j,n},$ $%
\I
(\cdot)$ denotes the indicator function, and the relevant inverses are assumed
to exist. Defining $\mathbf{\hat{\Gamma}}_{n}=\sum_{i=1}^{n}\mathbf{\hat{v}%
}_{i,n}\mathbf{\hat{v}}_{i,n}^{\prime}/n,$ the objective is to establish a
valid Gaussian distributional approximation of the finite sample distribution
of the OLS estimator $\mathbf{\hat{\beta}}_{n}$, and then find an estimator
$\mathbf{\hat{\Sigma}}_{n}$ of the variance of $\sum_{i=1}^{n}\mathbf{\hat{v}%
}_{i,n}u_{i,n}/\sqrt{n}$ such that%
\begin{equation}
\mathbf{\hat{\Omega}}_{n}^{-1/2}\sqrt{n}(\mathbf{\hat{\beta}}_{n}%
-\mathbf{\beta})\rightarrow_{d}\mathcal{N}(\mathbf{0},\mathbf{I}%
),\qquad\mathbf{\hat{\Omega}}_{n}=\mathbf{\hat{\Gamma}}_{n}^{-1}%
\mathbf{\hat{\Sigma}}_{n}\mathbf{\hat{\Gamma}}_{n}^{-1}%
,\label{eq:ANStudentizedEstimator}%
\end{equation}
in which case asymptotic valid inference on $\mathbf{\beta}$ can be conducted
in the usual way by employing the distributional approximation $\mathbf{\hat
{\beta}}_{n}\overset{a}{\sim}\mathcal{N}(\mathbf{\beta},\mathbf{\hat{\Omega}%
}_{n}/n)$. Our assumptions below will ensure that $\mathbf{\hat{\beta}}_{n}$
remains $\sqrt{n}$-consistent because we show in the supplemental appendix
that $\mathbf{\hat{\Omega}}_{n}^{-1}=O_{p}(1)$ even when $K_{n}%
/n\not \rightarrow 0$.

Our first result, Theorem \ref{thm:1} below, gives sufficient conditions for a
valid Gaussian approximation of the distribution of the infeasible statistic
$\mathbf{\Omega}_{n}^{-1/2}\sqrt{n}(\mathbf{\hat{\beta}}_{n}-\mathbf{\beta})$,
where $\mathbf{\Omega}_{n}=\mathbf{\hat{\Gamma}}_{n}^{-1}\mathbf{\Sigma}%
_{n}\mathbf{\hat{\Gamma}}_{n}^{-1}$ and $\mathbf{\Sigma}_{n} $ denotes the
variance of $\sum_{i=1}^{n}\mathbf{\hat{v}}_{i,n}u_{i,n}/\sqrt{n}$, even when
possibly $K_{n}/n\not \rightarrow 0$ and the linear regression model exhibits
conditional heteroskedasticity. This result, in turn, gives the basic
ingredient for discussing valid variance estimation in high-dimensional linear
regression models. Defining $\hat{u}_{i,n}=\sum_{j=1}^{n}M_{ij,n}%
(y_{j,n}-\mathbf{\hat{\beta}}_{n}^{\prime}\mathbf{x}_{j,n}),$ standard choices
of $\mathbf{\hat{\Sigma}}_{n}$ in the fixed-$K_{n}$ case include the
homoskedasticity-only estimator%
\[
\mathbf{\hat{\Sigma}}_{n}^{\mathtt{HO}}=\hat{\sigma}_{n}^{2}\mathbf{\hat
{\Gamma}}_{n},\qquad\hat{\sigma}_{n}^{2}=\frac{1}{n-d-K_{n}}\sum_{i=1}^{n}%
\hat{u}_{i,n}^{2},
\]
and the Eicker-White-type estimator%
\[
\mathbf{\hat{\Sigma}}_{n}^{\mathtt{EW}}=\frac{1}{n}\sum_{i=1}^{n}%
\mathbf{\hat{v}}_{i,n}\mathbf{\hat{v}}_{i,n}^{\prime}\hat{u}_{i,n}^{2}.
\]

Perhaps not too surprisingly, in Theorem \ref{thm:2} below, we find that
consistency of $\mathbf{\hat{\Sigma}}_{n}^{\mathtt{HO}}$ under
homoskedasticity holds quite generally even for models with many covariates.
In contrast, construction of a heteroskedasticity-robust estimator of
$\mathbf{\Sigma}_{n}$ is more challenging, as it turns out that consistency of
$\mathbf{\hat{\Sigma}}_{n}^{\mathtt{EW}}$ generally requires $K_{n}$ to be a
vanishing fraction of $n$.

To fix ideas, suppose $(y_{i,n},\mathbf{x}_{i,n}^{\prime},\mathbf{w}%
_{i,n}^{\prime})$ are i.i.d. over $i.$ It turns out that, under certain
regularity conditions,%
\[
\mathbf{\hat{\Sigma}}_{n}^{\mathtt{EW}}=\frac{1}{n}\sum_{i=1}^{n}\sum
_{j=1}^{n}M_{ij,n}^{2}\mathbf{\hat{v}}_{i,n}\mathbf{\hat{v}}_{i,n}^{\prime
}\mathbb{E}[u_{j,n}^{2}|\mathbf{x}_{j,n},\mathbf{w}_{j,n}]+o_{p}(1),
\]
whereas a requirement for (\ref{eq:ANStudentizedEstimator}) to hold is that
the estimator $\mathbf{\hat{\Sigma}}_{n}$ satisfies%
\begin{equation}
\mathbf{\hat{\Sigma}}_{n}=\frac{1}{n}\sum_{i=1}^{n}\mathbf{\hat{v}}%
_{i,n}\mathbf{\hat{v}}_{i,n}^{\prime}\mathbb{E}[u_{i,n}^{2}|\mathbf{x}%
_{i,n},\mathbf{w}_{i,n}]+o_{p}(1).\label{eq:AsyVar-Sigma}%
\end{equation}
The difference between the leading terms in the expansions is non-negligible
in general unless $K_{n}/n\rightarrow0.$ In recognition of this problem with
$\mathbf{\hat{\Sigma}}_{n}^{\mathtt{EW}},$ we study the more general class of
estimators of the form%
\[
\mathbf{\hat{\Sigma}}_{n}(\mathbf{\kappa}_{n})=\frac{1}{n}\sum_{i=1}^{n}%
\sum_{j=1}^{n}\kappa_{ij,n}\mathbf{\hat{v}}_{i,n}\mathbf{\hat{v}}%
_{i,n}^{\prime}\hat{u}_{j,n}^{2},
\]
where $\kappa_{ij,n}$ denotes element $\left(  i,j\right)  $ of a symmetric
matrix $\mathbf{\kappa}_{n}=\mathbf{\kappa}_{n}(\mathbf{w}_{1,n}%
,\ldots,\mathbf{w}_{n,n}).$ Estimators that can be written in this fashion
include $\mathbf{\hat{\Sigma}}_{n}^{\mathtt{EW}}$ (which corresponds to
$\mathbf{\kappa}_{n}=\mathbf{I}_{n}$) as well as variants of the so-called
HC$k$ estimators, $k\in\{1,2,3,4\}$, reviewed by \cite{Long-Ervin_2000_AS} and
\cite{MacKinnon_2012_BookCh}, among many others. To be specific, a natural
variant of HC$k$ is obtained by choosing $\kappa_{n}$ to be diagonal with
$\kappa_{ii,n}=\Upsilon_{i,n}M_{ii,n}^{-\xi_{i,n}}$, where $(\Upsilon
_{i,n},\xi_{i,n})=(1,0)$ for HC$0$ (and corresponding to $\mathbf{\hat{\Sigma
}}_{n}^{\mathtt{EW}}$), $(\Upsilon_{i,n},\xi_{i,n})=(n/(n-K_{n}),0)$ for
HC$1$, $(\Upsilon_{i,n},\xi_{i,n})=(1,1)$ for HC$2$, $(\Upsilon_{i,n}%
,\xi_{i,n})=(1,2)$ for HC$3$, and $(\Upsilon_{i,n},\xi_{i,n})=(1,\min
(4,nM_{ii,n}/K_{n}))$ for HC$4$. See Sections \ref{section:HCk} for more details.

In Theorem \ref{thm:3} below, we show that all of the HC$k$-type estimators,
which correspond to a diagonal choice of $\mathbf{\kappa}_{n}$, have the
shortcoming that they do not satisfy (\ref{eq:AsyVar-Sigma}) when
$K_{n}/n\nrightarrow0$. On the other hand, it turns out that a certain
non-diagonal choice of $\mathbf{\kappa}_{n}$ makes it possible to satisfy
(\ref{eq:AsyVar-Sigma}) even if $K_{n}$ is a non-vanishing fraction of $n$. To
be specific, it turns out that (under regularity conditions and) under mild
conditions under the weights $\kappa_{ij,n},$ $\mathbf{\hat{\Sigma}}%
_{n}(\mathbf{\kappa}_{n})$ satisfies%
\begin{equation}
\mathbf{\hat{\Sigma}}_{n}(\mathbf{\kappa}_{n})=\frac{1}{n}\sum_{i=1}^{n}%
\sum_{j=1}^{n}\sum_{k=1}^{n}\kappa_{ik,n}M_{kj,n}^{2}\mathbf{\hat{v}}%
_{i,n}\mathbf{\hat{v}}_{i,n}^{\prime}\mathbb{E}[u_{j,n}^{2}|\mathbf{x}%
_{j,n},\mathbf{w}_{j,n}]+o_{p}(1),\label{eq:ApproxSigmaHat}%
\end{equation}
suggesting that (\ref{eq:AsyVar-Sigma}) holds with $\mathbf{\hat{\Sigma}}%
_{n}=\mathbf{\hat{\Sigma}}_{n}(\mathbf{\kappa}_{n})$ provided $\kappa_{n}$ is
chosen in such a way that%
\begin{equation}
\sum_{k=1}^{n}\kappa_{ik,n}M_{kj,n}^{2}=%
\I
(i=j),\text{\qquad}1\leq i,j\leq n.\label{eq:kappasystem}%
\end{equation}
Accordingly, we define%
\[
\mathbf{\hat{\Sigma}}_{n}^{\mathtt{HC}}=\mathbf{\hat{\Sigma}}_{n}%
(\mathbf{\kappa}_{n}^{\mathtt{HC}})=\frac{1}{n}\sum_{i=1}^{n}\sum_{j=1}%
^{n}\kappa_{ij,n}^{\mathtt{HC}}\mathbf{\hat{v}}_{i,n}\mathbf{\hat{v}}%
_{i,n}^{\prime}\hat{u}_{j,n}^{2},
\]
where, with $\mathbf{M}_{n}$ denoting the matrix with element $(i,j)$ given by
$M_{ij,n}$ and $\odot$ denoting the Hadamard product,%
\[
\mathbf{\kappa}_{n}^{\mathtt{HC}}=\left(
\begin{array}
[c]{ccc}%
\kappa_{11,n}^{\mathtt{HC}} & \cdots & \kappa_{1n,n}^{\mathtt{HC}}\\
\vdots & \ddots & \vdots\\
\kappa_{n1,n}^{\mathtt{HC}} & \cdots & \kappa_{nn,n}^{\mathtt{HC}}%
\end{array}
\right)  =\left(
\begin{array}
[c]{ccc}%
M_{11,n}^{2} & \cdots & M_{1n,n}^{2}\\
\vdots & \ddots & \vdots\\
M_{n1,n}^{2} & \cdots & M_{nn,n}^{2}%
\end{array}
\right)  ^{-1}=(\mathbf{M}_{n}\odot\mathbf{M}_{n})^{-1}.
\]
The estimator $\mathbf{\hat{\Sigma}}_{n}^{\mathtt{HC}}$ is well defined
whenever $\mathbf{M}_{n}\odot\mathbf{M}_{n}$ is invertible, a simple
sufficient condition for which is that $\mathcal{M}_{n}<1/2,$ where%
\[
\mathcal{M}_{n}=1-\min_{1\leq i\leq n}M_{ii,n}.
\]
The fact that $\mathcal{M}_{n}<1/2$ implies invertibility of $\mathbf{M}%
_{n}\odot\mathbf{M}_{n}$ is a consequence of the Gershgorin circle theorem.
For details, see Section 3 in the supplemental appendix. More importantly, a
slight strengthening of the condition $\mathcal{M}_{n}<1/2$ will be shown to
be sufficient for (\ref{eq:ANStudentizedEstimator}) and (\ref{eq:AsyVar-Sigma}%
) to hold with $\mathbf{\hat{\Sigma}}_{n}=\mathbf{\hat{\Sigma}}_{n}%
^{\mathtt{HC}}$. Our final result, Theorem \ref{thm:4} below, formalizes this
finding (see also the supplemental appendix for further intuition underlying
this result).

The key intuition underlying our variance estimation result is that, even
though each conditional variance $\mathbb{E}[u_{i,n}^{2}|\mathbf{x}%
_{i,n},\mathbf{w}_{i,n}]$ cannot be well estimated due to the curse of
dimensionality, an averaged version such as the leading term in
(\ref{eq:AsyVar-Sigma}) can be estimated consistently. Thus, taking
$\widehat{\mathbb{E}}[u_{i,n}^{2}|\mathbf{x}_{i,n},\mathbf{w}_{i,n}%
]=\sum_{k=1}^{n}\kappa_{ik,n}\hat{u}_{k,n}^{2}$ as an estimator of
$\mathbb{E}[u_{i,n}^{2}|\mathbf{x}_{i,n},\mathbf{w}_{i,n}]$, plugging into the
leading term in (\ref{eq:AsyVar-Sigma}), and computing conditional
expectations, we obtain the leading term in (\ref{eq:ApproxSigmaHat}). To make
this leading term equal to the desired target $\sum_{i=1}^{n}\mathbf{\hat{v}%
}_{i,n}\mathbf{\hat{v}}_{i,n}^{\prime}\mathbb{E}[u_{i,n}^{2}|\mathbf{x}%
_{i,n},\mathbf{w}_{i,n}]$, it is natural to require
\[
\sum_{j=1}^{n}\sum_{k=1}^{n}\kappa_{ik,n}M_{kj,n}^{2}\mathbb{E}[u_{j,n}%
^{2}|\mathbf{x}_{j,n},\mathbf{w}_{j,n}]=\mathbb{E}[u_{i,n}^{2}|\mathbf{x}%
_{i,n},\mathbf{w}_{i,n}]\text{\qquad}1\leq i\leq n.
\]
Since $\mathbb{E}[u_{i,n}^{2}|\mathbf{x}_{i,n},\mathbf{w}_{i,n}]$ are unknown,
our variance estimator solves (\ref{eq:kappasystem}), which generates enough
equations to solve for all $n(n-1)/2$ possibly distinct elements in
$\mathbf{\kappa}_{n}^{\mathtt{HC}}$.

\begin{description}
\item[Remark 1.] $\mathbf{\hat{\Sigma}}_{n}^{\mathtt{HC}}=n^{-1}\sum_{i=1}%
^{n}\mathbf{\hat{v}}_{i,n}\mathbf{\hat{v}}_{i,n}^{\prime}\tilde{u}_{i,n}^{2}$
with $\tilde{u}_{i,n}^{2}=\sum_{j=1}^{n}\kappa_{ij,n}^{\mathtt{HC}}\hat
{u}_{j,n}^{2}$, and therefore $\tilde{u}_{i,n}^{2}$ can be interpreted as a
bias-corrected \textquotedblleft estimator\textquotedblright\ of (the
conditional expectation of) $u_{i,n}^{2} $.
\end{description}

\section{Setup\label{section:framework}}

This section introduces a general framework encompassing several special cases
of linear-in-parameters regression models of the form (\ref{eq:OLSmodel}). We
first present generic high-level assumptions, and then discuss their
implications as well as some easier to verify sufficient conditions. Finally,
to close this setup section, we briefly discuss three motivating leading
examples: linear regression models with increasing dimension, muti-way fixed
effect linear models, and semiparametric semi-linear regression. Technical
details and related results for these examples are given in the supplemental appendix.

\subsection{Framework}

Suppose $\{(y_{i,n},\mathbf{x}_{i,n}^{\prime},\mathbf{w}_{i,n}^{\prime}):1\leq
i\leq n\}$ is generated by (\ref{eq:OLSmodel}). Let $\Vert\cdot\Vert$ denote
the Euclidean norm, set $\mathcal{X}_{n}=(\mathbf{x}_{1,n},\ldots
,\mathbf{x}_{n,n})$, and for a collection $\mathcal{W}_{n}$ of random
variables satisfying $\mathbb{E}[\mathbf{w}_{i,n}|\mathcal{W}_{n}%
]=\mathbf{w}_{i,n}$, define the constants%
\begin{align*}
\varrho_{n}  & =\frac{1}{n}\sum_{i=1}^{n}\mathbb{E}[R_{i,n}^{2}],\qquad
R_{i,n}=\mathbb{E}[u_{i,n}|\mathcal{X}_{n},\mathcal{W}_{n}],\\
\rho_{n}  & =\frac{1}{n}\sum_{i=1}^{n}\mathbb{E}[r_{i,n}^{2}],\qquad
r_{i,n}=\mathbb{E}[u_{i,n}|\mathcal{W}_{n}],\\
\chi_{n}  & =\frac{1}{n}\sum_{i=1}^{n}\mathbb{E}[\Vert\mathbf{Q}_{i,n}%
\Vert^{2}],\qquad\mathbf{Q}_{i,n}=\mathbb{E}[\mathbf{v}_{i,n}|\mathcal{W}%
_{n}],
\end{align*}
where $\mathbf{v}_{i,n}=\mathbf{x}_{i,n}-(\sum_{j=1}^{n}\mathbb{E}%
[\mathbf{x}_{j,n}\mathbf{w}_{j,n}^{\prime}])(\sum_{j=1}^{n}\mathbb{E}%
[\mathbf{w}_{j,n}\mathbf{w}_{j,n}^{\prime}])^{-1}\mathbf{w}_{i,n}$ is the
population counterpart of $\mathbf{\hat{v}}_{i,n}.\ $Also, define%
\[
\mathcal{C}_{n}=\max_{1\leq i\leq n}\{\mathbb{E}[U_{i,n}^{4}|\mathcal{X}%
_{n},\mathcal{W}_{n}]+\mathbb{E}[\Vert\mathbf{V}_{i,n}\Vert^{4}|\mathcal{W}%
_{n}]+1/\mathbb{E}[U_{i,n}^{2}|\mathcal{X}_{n},\mathcal{W}_{n}]\}+1/\lambda
_{\min}(\mathbb{E}[\mathbf{\tilde{\Gamma}}_{n}|\mathcal{W}_{n}])\},
\]
where $U_{i,n}=y_{i,n}-\mathbb{E}[y_{i,n}|\mathcal{X}_{n},\mathcal{W}_{n}],$
$\mathbf{V}_{i,n}=\mathbf{x}_{i,n}-\mathbb{E}[\mathbf{x}_{i,n}|\mathcal{W}%
_{n}],$ $\mathbf{\tilde{\Gamma}}_{n}=\sum_{i=1}^{n}\mathbf{\tilde{V}}%
_{i,n}\mathbf{\tilde{V}}_{i,n}^{\prime}/n,$ and $\mathbf{\tilde{V}}_{i,n}%
=\sum_{j=1}^{n}M_{ij,n}\mathbf{V}_{j,n}.$

We impose the following three high-level conditions. Let $\lambda_{\min}%
(\cdot)$ denote the minimum eigenvalue of its argument, and $\overline{\lim
}_{n\rightarrow\infty}a_{n}=\lim\sup_{n\rightarrow\infty}a_{n}$ for any
sequence $a_{n}$.

\begin{assumption}
[Sampling]\label{ass:sampling}$\mathbb{C}[U_{i,n},U_{j,n}|\mathcal{X}%
_{n},\mathcal{W}_{n}]=0$ for $i\neq j$ and $\max_{1\leq i\leq N_{n}%
}\#\mathcal{T}_{i,n}=O(1), $ where $\#\mathcal{T}_{i,n}$ is the cardinality of
$\mathcal{T}_{i,n}$ and where $\{\mathcal{T}_{i,n}:1\leq i\leq N_{n}\}$ is a
partition of $\{1,\ldots,n\}$ such that $\{(U_{t,n},V_{t,n}):t\in
\mathcal{T}_{i,n}\}$ are independent over $i$ conditional on $\mathcal{W}%
_{n}.$
\end{assumption}

\begin{assumption}
[Design]\label{ass:design}$\mathbb{P}[\lambda_{\min}(\sum_{i=1}^{n}%
\mathbf{w}_{i,n}\mathbf{w}_{i,n}^{\prime})>0]\rightarrow1,$ $\overline{\lim
}_{n\rightarrow\infty}K_{n}/n<1,$ and $\mathcal{C}_{n}=O_{p}(1).$
\end{assumption}

\begin{assumption}
[Approximations]\label{ass:approximations}$\chi_{n}=O(1),$ $\varrho
_{n}+n(\varrho_{n}-\rho_{n})+n\chi_{n}\varrho_{n}=o(1),$ and $\max_{1\leq
i\leq n}\Vert\mathbf{\hat{v}}_{i,n}\Vert/\sqrt{n}=o_{p}(1).$
\end{assumption}

\subsection{Discussion of Assumptions}

Assumptions \ref{ass:sampling}--\ref{ass:approximations} are meant to be
high-level and general, allowing for different linear-in-parameters regression
models. We now discuss the main restrictions imposed by these assumptions. We
further illustrate them in the following subsection using more specific examples.

\subsubsection{Assumption \ref{ass:sampling}}

This assumption concerns the sampling properties of the observed data. It
generalizes classical i.i.d. sampling by allowing for groups or
\textquotedblleft clusters\textquotedblright\ of finite but possibly
heterogeneous size with arbitrary intra-group dependence, which is very common
in the context of fixed effects linear regression models. As currently stated,
this assumption does not allow for dependence in the error terms across units,
and therefore excludes clustered, spacial or time series dependence in the
sample. We conjecture our main results extend to the latter cases, though here
we focus on i.n.i.d. (conditionally)\ heteroskedastic models only, and hence
relegate the extension to errors exhibiting clustered, spacial or time series
dependence for future work. Assumption \ref{ass:sampling} reduces to classical
i.i.d. sampling when $N_{n}=n$, $\mathcal{T}_{i,n}=\{i\}$ [implying
$\max_{1\leq i\leq N_{n}}\#\mathcal{T}_{i,n}=1$], and all observations have
the same distribution.

\subsubsection{Assumption \ref{ass:design}}

This assumption concerns basic design features of the linear regression model.
The first two restrictions are mild and reflect the main goal of this paper,
that is, analyzing linear regression models with many nuisance covariates
$\mathbf{w}_{i,n}$. In practice, the first restriction regarding the minimum
eigenvalue of the design matrix $\sum_{i=1}^{n}\mathbf{w}_{i,n}\mathbf{w}%
_{i,n}^{\prime}$ is always imposed by removing redundant (i.e., linearly
dependent) covariates; from a theoretical perspective this condition requires
either restrictions on the distributional relationship of such covariates or
some form of trimming leading to selection of included covariates (e.g., most
software packages remove covariates leading to \textquotedblleft
too\textquotedblright\ small eigenvalues of the design matrix by means of some
hard-thresholding rule).

On the other hand, the last condition, $\mathcal{C}_{n}=O_{p}(1),$ may be
restrictive in some settings: for example, if the covariates have unbounded
support (e.g., they are normally distributed) and heteroskedasticity is
unbounded (e.g., unbounded multiplicative heteroskedasticity), then the
assumption may fail. Simple sufficient conditions for $\mathcal{C}_{n}%
=O_{p}(1)$ can be formulated when the covariates have compact support, or the
heteroskedasticity is multiplicative and bounded, because in these cases it is
easy to bound the conditional moments of the error terms. It would be useful
to know whether the condition $\mathcal{C}_{n}=O_{p}(1)$ can be relaxed to a
version involving only unconditional moments, though we conjecture this
weaker\ assumption will require a different method of proof (see the
supplemental appendix for details).

\subsubsection{Assumption \ref{ass:approximations}}

This assumption requires two basic approximations to hold. First, concerning
bias, conditions on $\varrho_{n}$ are related to the approximation quality of
the linear-in-parameters model (\ref{eq:OLSmodel}) for the \textquotedblleft
long\textquotedblright\ conditional expectation $\mathbb{E}[y_{i,n}%
|\mathcal{X}_{n},\mathcal{W}_{n}]$. Similarly, conditions on $\rho_{n}$ and
$\chi_{n}$ are related to linear-in-parameters approximations for the
\textquotedblleft short\textquotedblright\ conditional expectations
$\mathbb{E}[y_{i,n}|\mathcal{W}_{n}]$ and $\mathbb{E}[\mathbf{x}%
_{i,n}|\mathcal{W}_{n}]$, respectively. All these approximations are measured
in terms of population mean square error, and are at the heart of empirical
work employing linear-in-parameters regression models. Depending on the model
of interest, different sufficient conditions can be given for these
assumptions. Here we briefly mention the most simple one:\ (\emph{a}) if
$\mathbb{E}[u_{i,n}|\mathcal{X}_{n},\mathcal{W}_{n}]=0$ for all $i$ and $n$,
which can be interpreted as exogeneity (e.g., no misspecification bias), then
$0=\rho_{n}=n(\varrho_{n}-\rho_{n})+n\chi_{n}\varrho_{n}$ for all $n$; and
(\emph{b})\ if $\mathbb{E}[\Vert\mathbf{x}_{i,n}\Vert^{2}]<\infty$ for all $i$
and $n$, then $\chi_{n}=O(1)$. Other sufficient conditions are discussed below.

Second, the high-level condition $\max_{1\leq i\leq n}\Vert\mathbf{\hat{v}%
}_{i,n}\Vert/\sqrt{n}=o_{p}(1)$ restricts the distributional relationship
between the finite dimensional covariate of interest $\mathbf{x}_{i,n}$ and
the high-dimensional nuisance covariate $\mathbf{w}_{i,n}$. This condition can
be interpreted as a negligibility condition and thus comes close to minimal
for the central limit theorem to hold. At the present level of generality it
seems difficult to formulate primitive sufficient conditions for this
restriction that cover all cases of interest, but for completeness we mention
that under mild moment conditions it suffices to require that one of the
following conditions hold (see Lemma SA-7 in the supplemental appendix for
details and weaker conditions):\newline\textbf{(i)} $\mathcal{M}_{n}%
=o_{p}(1),$ or\newline\textbf{(ii)} $\chi_{n}=o(1),$ or \newline\textbf{(iii)}
$\max_{1\leq i\leq n}\sum_{j=1}^{n}%
\I
(M_{ij,n}\neq0)=o_{p}(n^{1/3})$.

Each of these conditions is interpretable. First, $\mathcal{M}_{n}\geq
K_{n}/n$ because $\sum_{i=1}^{n}M_{ii,n}=n-K_{n}$ and a necessary condition
for \textbf{(i)} is therefore that $K_{n}/n\rightarrow0.$ Conversely, because%
\[
\mathcal{M}_{n}\leq\frac{K_{n}}{n}\frac{1-\min_{1\leq i\leq n}M_{ii,n}}%
{1-\max_{1\leq i\leq n}M_{ii,n}},
\]
the condition $K_{n}/n\rightarrow0$ is sufficient for \textbf{(i)} whenever
the design is \textquotedblleft approximately balanced\textquotedblright\ in
the sense that $(1-\min_{1\leq i\leq n}M_{ii,n})/(1-\max_{1\leq i\leq
n}M_{ii,n})=O_{p}(1).$ In other words, \textbf{(i)} requires and effectively
covers the case where it is assumed that $K_{n}$ is a vanishing fraction of
$n.$ In contrast, conditions \textbf{(ii)} and \textbf{(iii)} can hold also
when $K_{n}$ is a non-vanishing fraction of $n,$ which is the case of primary
interest in this paper.

Because \textbf{(ii)} is a requirement on the accuracy of the approximation
$\mathbb{E}[\mathbf{x}_{i,n}|\mathbf{w}_{i,n}]\approx\mathbf{\delta}%
_{n}^{\prime}\mathbf{w}_{i,n}$ with $\mathbf{\delta}_{n}=\mathbb{E}%
[\mathbf{w}_{i,n}\mathbf{w}_{i,n}^{\prime}]^{-1}\mathbb{E}[\mathbf{w}%
_{i,n}\mathbf{x}_{i,n}^{\prime}]$, primitive conditions for it are available
when, for example, the elements of $\mathbf{w}_{i,n}$ are approximating
functions. Indeed, in such cases one typically has $\chi_{n}=O(K_{n}^{-\alpha
})$ for some $\alpha>0$, so condition \textbf{(ii)} not only accommodates
$K_{n}/n\nrightarrow0,$ but actually places no upper bound on the magnitude of
$K_{n}$ in important special cases. This condition also holds when
$\mathbf{w}_{i,n}$ are dummy variables or discrete covariates, as we discuss
in more detail below.

Finally, condition \textbf{(iii)}, and its underlying higher-level condition
described in the supplemental appendix, is useful to handle cases where
$\mathbf{w}_{i,n}$ cannot be interpreted as approximating functions, but
rather just many different covariates included in the linear model
specification. This condition is a \textquotedblleft
sparsity\textquotedblright\ condition on the projection matrix $\mathbf{M}%
_{n}$, which allows for $K_{n}/n\nrightarrow0$. The condition is easy to
verify in certain cases, including those where \textquotedblleft locally
bounded\textquotedblright\ approximating functions or fixed effects are used
(see below for concrete examples).

\subsection{Motivating Examples}

We briefly mention three motivating examples of linear-in-parameter regression
models covered by our results. All technical details are given in the
supplemental appendix.

\subsubsection{Linear Regression Model with Increasing Dimension}

This leading example has a long tradition in statistics and econometrics. The
model takes (\ref{eq:OLSmodel}) as the data generating process (DGP),
typically with i.i.d. data and the exogeneity condition $\mathbb{E}%
[u_{i,n}|\mathbf{x}_{i,n},\mathbf{w}_{i,n}]=0$. However, our assumptions only
require $n\mathbb{E}[(\mathbb{E}[u_{i,n}|\mathbf{x}_{i,n},\mathbf{w}%
_{i,n}])^{2}]=o(1)$, and hence (\ref{eq:OLSmodel}) can be interpreted as a
linear-in-parameters mean-square approximation to the unknown conditional
expectation $\mathbb{E}[y_{i,n}|\mathbf{x}_{i,n},\mathbf{w}_{i,n}]$. Either
way, $\mathbf{\hat{\beta}}_{n}$ is the standard OLS estimator.

Setting $\mathcal{W}_{n}=(\mathbf{w}_{1,n},\ldots,\mathbf{w}_{n,n}),$
$N_{n}=n,$ $\mathcal{T}_{i,n}=\{i\}$ and $\max_{1\leq i\leq N_{n}%
}\#\mathcal{T}_{i,n}=1$, Assumptions \ref{ass:sampling}--\ref{ass:design} are
standard, while Assumption \ref{ass:approximations} is satisfied provided that
$\mathbb{E}[\Vert\mathbf{x}_{i,n}\Vert^{2}]<\infty$ [implying $\chi_{n}%
=O(1)$], $n\mathbb{E}[(\mathbb{E}[u_{i,n}|\mathbf{x}_{i,n},\mathbf{w}%
_{i,n}])^{2}]=o(1)$ [implying $n(\varrho_{n}-\rho_{n})+n\chi_{n}\varrho
_{n}=o(1)$], and $\max_{1\leq i\leq n}\Vert\mathbf{\hat{v}}_{i,n}\Vert
/\sqrt{n}=o_{p}(1)$. Primitive sufficient conditions for the latter
negligibility condition can be given as discussed above. For example, under
regularity conditions, $\chi_{n}=o(1)$ if either (\emph{a})$\ \mathbb{E}%
[\mathbf{x}_{i,n}|\mathbf{w}_{i,n}]=\mathbf{\delta}^{\prime}\mathbf{w}_{i,n}
$, (\emph{b})$\ $the nuisance covariates are discrete and a saturated dummy
variables model is used, or (\emph{c}) $\mathbf{w}_{i,n}$ are constructed
using sieve functions. Alternatively, $\max_{1\leq i\leq n}\sum_{j=1}^{n}%
\I
(M_{ij,n}\neq0)=o_{p}(n^{1/3})$ is satisfied provided the distribution of the
nuisance covariates $\mathbf{w}_{i,n}$ generates a projection matrix
$\mathbf{M}_{n}$ that is approximately a band matrix (see below for concrete
examples). Precise regularity conditions for this example are given in Section
4.1 of the supplemental appendix.

\subsubsection{Fixed Effects Panel Data Regression Model}

A second class of examples covered by our results are linear panel data models
with multi-way fixed effects and related models such as those encountered in
networks, spillovers or social interactions settings. A common feature in
these examples is the presence of possibly many dummy variables in
$\mathbf{w}_{i,n}$, capturing unobserved heterogeneity or other unobserved
effects across units (e.g., network link or spillover effect). In many
applications the number of distinct dummy-type variables is large because
researchers often include multi-group indicators, interactions thereof, and
similar regressors obtained from factor variables. In these complicated models
the nuisance covariates need to be estimated explicitly, even in simple linear
regression problems, because it is not possible to difference out the
multi-way indicator variables for estimation and inference.

\cite{Stock-Watson_2008_ECMA} consider heteroskedasticity-robust inference for
the one-way fixed effect panel data regression model%
\begin{equation}
Y_{it}=\alpha_{i}+\mathbf{\beta}^{\prime}\mathbf{X}_{it}+U_{it},\text{\qquad
}i=1,\ldots,N,\text{\qquad}t=1,\ldots,T,\label{eq:PanelDataModel}%
\end{equation}
where $\alpha_{i}\in\mathbb{R}$ is an individual-specific intercept,
$\mathbf{X}_{it}$ is a regressor of dimension $d$, and $U_{it}$ is an scalar
error term, and the following assumptions are satisfied. To map this model
into our framework, suppose that $\{(U_{i1},\ldots,U_{iT},\mathbf{X}%
_{i1}^{\prime}\ldots,\mathbf{X}_{iT}^{\prime}):1\leq i\leq n\}$ are
independent over $i$, $\mathbb{E}[U_{it}|\mathbf{X}_{i1}\ldots,\mathbf{X}%
_{iT}]=0$, and $\mathbb{E}[U_{it}U_{is}|X_{i1}\ldots,X_{iT}]=0$ for $t\neq s.$
Then, setting $n=NT$, $K_{n}=N$, $\mathbf{\gamma}_{n}=(\alpha_{1}%
,\ldots,\alpha_{N})^{\prime}$, and $(y_{(i-1)T+t,n},\mathbf{x}_{(i-1)T+t,n}%
^{\prime},u_{(i-1)T+t,n},\mathbf{w}_{(i-1)T+t,n}^{\prime})=(Y_{it}%
,\mathbf{X}_{it}^{\prime},U_{it},\mathbf{e}_{i,N}^{\prime})$, $1\leq i\leq N$
and $1\leq t\leq T$, where $\mathbf{e}_{i,N}\in\mathbb{R}^{N} $ is the $i$-th
unit vector of dimension $N$, the model (\ref{eq:PanelDataModel}) is also of
the form (\ref{eq:OLSmodel}) and $\mathbf{\hat{\beta}}_{n}$ is the fixed
effects estimator of $\mathbf{\beta}$. In general, this model does not satisfy
an i.i.d. assumption, but Assumption \ref{ass:sampling} enables us to employ
results for independent random variables when developing asymptotics. In
particular, unlike \cite{Stock-Watson_2008_ECMA}, we do not require
$(U_{i1},\ldots,U_{iT},\mathbf{X}_{i1}^{\prime}\ldots,\mathbf{X}_{iT}^{\prime
})$ to be i.i.d. over $i$, nor we require any kind of stationarity on the part
of $(U_{it},\mathbf{X}_{it}^{\prime})$. The amount of variance heterogeneity
permitted is quite large, since we basically only require $\mathbb{V}%
[Y_{it}|\mathbf{X}_{i1},\ldots,\mathbf{X}_{iT}]=\mathbb{E}[U_{it}%
^{2}|\mathbf{X}_{i1},\ldots,\mathbf{X}_{iT}]$ to be bounded and bounded away
from zero. (On the other hand, serial correlation is assumed away because our
assumptions imply that $\mathbb{C}[Y_{it},Y_{is}|\mathbf{X}_{i1}%
,\ldots,\mathbf{X}_{iT}]=0$ for $t\neq s.$) In other respects this model is in
fact more tractable than the previous models due to the special nature of the
covariates $\mathbf{w}_{i,n} $, that is, a dummy variable for each unit
$i=1,\ldots,N$.

In this one-way fixed effects example, $K_{n}/n=1/T$ and therefore a
high-dimensional model corresponds to a short panel model: $\max_{1\leq i\leq
n}\sum_{j=1}^{n}%
\I
(M_{ij,n}\neq0)=T$ and hence the negligibility condition holds easily. If
$T\geq2$, our asymptotic Gaussian approximation for the distribution of the
least-squares estimator $\mathbf{\hat{\beta}}_{n}$ is valid (see Theorem
\ref{thm:1}), despite the coefficients $\mathbf{\gamma}_{n}$ not being
consistently estimated. On the other hand, consistency of our generic variance
estimator requires $T\geq3$ [implying $K_{n}/n<1/2$]; see Theorems \ref{thm:3}
and \ref{thm:4}. Further details are given in Section 4.2 of the supplemental
appendix, where we also discuss a case-specific consistent variance estimator
when $T=2$.

Our generic results go beyond one-way fixed effect linear regression models,
as they can be used to obtain valid inference in other contexts where
multi-way fixed effects or similar discrete regressors are included. For a
second concrete example, consider the recent work of \cite[and references
therein]{Verdier_2017_wp} in the context of linear models with two-way
unobserved heterogeneity and sparsely matched data. This model is isomorphic
to a network model, where students and teacher (or workers and firms, for
another example) are \textquotedblleft matched\textquotedblright\ or
\textquotedblleft connected\textquotedblright\ over time, but potential
unobserved heterogeneity at both levels is a concern. In this setting, under
random sampling, \cite{Verdier_2017_wp} offers primitive conditions for our
high-level assumptions when two-way fixed effect models are used for
estimation and inference. In particular, using a clever Markov chain argument
(see his Lemma 1), he is able to provide different restriction on $T $ and the
number of matches in the network to ensure consistent variance estimation
using the methods developed in this paper. To give one concrete example, he
finds that if $T\geq5$ and for any pair of teachers (firms), the number of
students (workers) assigned to both teachers (firms) in the pair is either
zero or greater than three, then our key high-level condition in Theorem
\ref{thm:4} below is verified.

\subsubsection{Semiparametric Partially Linear Model}

Another model covered by our results is the partially linear model%
\begin{equation}
y_{i}=\mathbf{\beta}^{\prime}\mathbf{x}_{i}+g(\mathbf{z}_{i})+\varepsilon
_{i},\qquad i=1,\ldots,n,\label{eq:PartiallyLinearModel}%
\end{equation}
where $\mathbf{x}_{i}$ and $\mathbf{z}_{i}$ are explanatory variables,
$\varepsilon_{i}$ is an error term satisfying $\mathbb{E}[\varepsilon
_{i}|\mathbf{x}_{i},\mathbf{z}_{i}]=0$, the function $g(\mathbf{z})$ is
unknown, and sampling is i.i.d. across $i$ is assumed. Suppose $\{\mathbf{p}%
^{k}(\mathbf{z}):k=1,2,\cdots\}$ are functions having the property that linear
combinations can approximate square-integrable functions of $\mathbf{z}$ well,
in which case $g(\mathbf{z}_{i})\approx\mathbf{\gamma}_{n}^{\prime}%
\mathbf{p}_{n}(\mathbf{z}_{i})$ for some $\mathbf{\gamma}_{n},$ where
$\mathbf{p}_{n}(\mathbf{z})=(\mathbf{p}^{1}(\mathbf{z}),\ldots,\mathbf{p}%
^{K_{n}}(\mathbf{z}))^{\prime}$. Defining $y_{i,n}=y_{i},$ $\mathbf{x}%
_{i,n}=\mathbf{x}_{i},$ $\mathbf{w}_{i,n}=\mathbf{p}_{n}(\mathbf{z}_{i}),$ and
$u_{i,n}=\varepsilon_{i}+g(\mathbf{z}_{i})-\mathbf{\gamma}_{n}^{\prime
}\mathbf{w}_{i,n},$ the model (\ref{eq:PartiallyLinearModel}) is of the form
(\ref{eq:OLSmodel}), and $\mathbf{\hat{\beta}}_{n}$ is the series estimator of
$\mathbf{\beta}$; see, e.g., \cite{Donald-Newey_1994_JMA} and
\cite{Cattaneo-Jansson-Newey_2017_ET} and references therein.

Constructing the basis $\mathbf{p}_{n}(\mathbf{z}_{i})$ in applications may
require using a large $K_{n}$, either because the underlying functions are not
smooth enough or because $\dim(\mathbf{z}_{i})$ is large. For example, if a
$\mathfrak{p}=3$ cubic polynomial expansion is used, also known as a power
series of order $3$, then $\dim(\mathbf{w}_{i})=(\mathfrak{p}+\dim
(\mathbf{z}_{i}))!/(\mathfrak{p}!\dim(\mathbf{z}_{i})!)=286$ when
$\dim(\mathbf{z}_{i})=10$, and therefore flexible estimation and inference
using the semi-linear model (\ref{eq:PartiallyLinearModel}) with a sample size
of $n=1,000$ gives $K_{n}/n=0.286$. For further technical details on
series-based methods see, e.g., \cite{Newey_1997_JoE},
\cite{Chen_2007_Handbook}, \cite{Cattaneo-Farrell_2013_JoE}, and
\cite{Belloni-Chernozhukov-Chetverikov-Kato_2015_JoE}, and references therein.
For another example, when the basis functions $\mathbf{p}_{n}(\mathbf{z})$ are
constructed using partitioning estimators, the OLS estimator of $\mathbf{\beta
}$ becomes a subclassification estimator, a method that has been proposed in
the literature on program evaluation and treatment effects; see, e.g.,
\cite{Cochran1968_Biometrics}, \cite{Rosenbaum-Rubin_1983_Biometrika},
\cite{Cattaneo-Farrell_2011_BookCh}, and references therein. When a
Partitioning estimator of order $0$ is used, the semi-linear model becomes a
one-way fixed effects linear regression model, where each dummy variable
corresponds to one (disjoint)\ partition on the support of $\mathbf{z}_{i}$;
in this case, $K_{n}$ is to the number of partitions or fixed effects included
in the estimation.

Our primitive regularity conditions for this example include%
\[
\varrho_{n}=\min_{\mathbf{\gamma}\in\mathbb{R}^{K_{n}}}\mathbb{E}%
[|g(\mathbf{z}_{i})-\mathbf{\gamma}^{\prime}\mathbf{p}_{n}(\mathbf{z}%
_{i})|^{2}]=o(1),\qquad\chi_{n}=\min_{\mathbf{\delta}\in\mathbb{R}%
^{K_{n}\times d}}\mathbb{E}[\Vert\mathbb{E}[\mathbf{x}_{i}|\mathbf{z}%
_{i}]-\mathbf{\delta}^{\prime}\mathbf{p}_{n}(\mathbf{z}_{i})\Vert^{2}]=O(1),
\]
$n\varrho_{n}\chi_{n}=o(1)$, and the negligibility condition $\max_{1\leq
i\leq n}\Vert\mathbf{\hat{v}}_{i,n}\Vert/\sqrt{n}=o_{p}(1)$. A key finding
implied by these regularity conditions is that we only require minimal
smoothness conditions on $g(\mathbf{z}_{i})$ and $\mathbb{E}[\mathbf{x}%
_{i}|\mathbf{z}_{i}]$. The negligibility condition is automatically satisfied
if $\chi_{n}=o(1)$, as discussed above, but in fact our results do not require
any approximation of $\mathbb{E}[\mathbf{x}_{i}|\mathbf{z}_{i}]$, as usually
assumed in the literature, provided a \textquotedblleft locally
supported\textquotedblright\ basis is used; i.e., any basis $\mathbf{p}%
_{n}(\mathbf{z})$ that generates an approximately band projection matrix
$\mathbf{M}_{n}$; examples of such basis include partitioning and spline
estimators. See Section 4.3 in the supplemental appendix for further
discussion and technical details.

\section{Results\label{section:results}}

This section presents our main theoretical results for inference in linear
regression models with many covariates and heteroskedasticity. Mathematical
proofs, and other technical results that may be of independent interest, are
given in the supplemental appendix.

\subsection{Asymptotic Normality}

As a means to the end of establishing (\ref{eq:ANStudentizedEstimator}), we
give an asymptotic normality result for $\mathbf{\hat{\beta}}_{n}$ which may
be of interest in its own right.

\begin{theorem}
\label{thm:1}Suppose Assumptions \ref{ass:sampling}--\ref{ass:approximations}
hold. Then,%
\begin{equation}
\mathbf{\Omega}_{n}^{-1/2}\sqrt{n}(\mathbf{\hat{\beta}}_{n}-\mathbf{\beta
})\rightarrow_{d}\mathcal{N}(\mathbf{0},\mathbf{I}),\text{\qquad
}\mathbf{\Omega}_{n}=\mathbf{\hat{\Gamma}}_{n}^{-1}\mathbf{\Sigma}%
_{n}\mathbf{\hat{\Gamma}}_{n}^{-1},\label{eq:ANStandarizedEstimator}%
\end{equation}
where $\mathbf{\Sigma}_{n}=\sum_{i=1}^{n}\mathbf{\hat{v}}_{i,n}\mathbf{\hat
{v}}_{i,n}^{\prime}\mathbb{E}[U_{i,n}^{2}|\mathcal{X}_{n},\mathcal{W}_{n}]/n.
$
\end{theorem}

In the literature on high-dimensional linear models, \cite{Mammen_1993_AoS}
obtains a similar asymptotic normality result as in Theorem \ref{thm:1} but
under the condition $K_{n}^{1+\delta}/n\rightarrow0$ for $\delta>0$ restricted
by certain moment condition on the covariates. In contrast, our result only
requires $\overline{\lim}_{n\rightarrow\infty}K_{n}/n<1,$ but imposes a
different restriction on the high-dimensional covariates (e.g., condition
\textbf{(i)}, \textbf{(ii)} or \textbf{(iii)} discussed previously) and
furthermore exploits the fact that the parameter of interest is given by the
first $d$ coordinates of the vector $(\mathbf{\beta}^{\prime},\mathbf{\gamma
}_{n}^{\prime})^{\prime}$ (i.e., in \cite{Mammen_1993_AoS} notation, it
considers the case $\mathbf{c}=(\mathbf{\iota}^{\prime},\mathbf{0}^{\prime
})^{\prime}$ with $\iota$ denoting a $d$-dimensional vector of ones and
$\mathbf{0}$ denoting a $K_{n}$-dimensional vector of zeros).

In isolation, the fact that Theorem \ref{thm:1} removes the requirement
$K_{n}/n\rightarrow0$ may seem like little more than a subtle technical
improvement over results currently available. It should be recognized,
however, that conducting inference turn out to be considerably harder when
$K_{n}/n\not \rightarrow 0.$ The latter is an important insight about
large-dimensional models that cannot be deduced from results obtained under
the assumption $K_{n}/n\rightarrow0,$ but can be obtained with the help of
Theorem \ref{thm:1}. In addition, it is worth mentioning that Theorem
\ref{thm:1} is a substantial improvement over \cite[Theorem 1]%
{Cattaneo-Jansson-Newey_2017_ET} because here it is not required that
$K_{n}\rightarrow\infty$ nor $\chi_{n}=o(1)$ ---a different method of proof is
also used. This improvement applies not only to the partially linear model
example, but more generally to linear models with many covariates, because
Theorem \ref{thm:1} applies to quite general form of nuisance covariate
$\mathbf{w}_{i,n}$ beyond specific approximating basis functions. In the
specific case of the partially linear model, this implies that we are able to
weaken smoothness assumptions (or the curse of dimensionality), otherwise
required to satisfy the condition $\chi_{n}=o(1)$.

\begin{description}
\item[Remark 2.] Theorem \ref{thm:1} concerns only distributional properties
of $\mathbf{\hat{\beta}}_{n}$. First, this theorem implies $\sqrt{n}%
$-consistency of $\mathbf{\hat{\beta}}_{n}$ because $\mathbf{\Omega}_{n}%
^{-1}=O_{p}(1)$ (see Lemmas SA-1 and SA-2 of the supplemental appendix).
Second, this theorem does require nor imply consistency of the
(implicit)\ least squares estimate of $\mathbf{\gamma}_{n}$, as in fact such a
result will not be true in most applications with many nuisance covariates
$\mathbf{w}_{n,i}$. For example, in a partially linear model
(\ref{eq:PartiallyLinearModel}) the approximating coefficients $\mathbf{\gamma
}_{n}$ will not be consistently estimated unless $K_{n}/n\rightarrow0$, or in
a one-way fixed effect panel data model (\ref{eq:PanelDataModel}) the
unit-specific coefficients in $\mathbf{\gamma}_{n}$ will not be consistently
estimated unless $K_{n}/n=1/T\rightarrow0$. Nevertheless, Theorem \ref{thm:1}
shows that $\mathbf{\hat{\beta}}_{n}$ can still be root-$n$ asymptotically
normal under fairly general conditions; this result is due to the intrinsic
linearity and additive separability of the model (\ref{eq:OLSmodel}).
\end{description}

\subsection{Variance Estimation}

Achieving (\ref{eq:ANStudentizedEstimator}), the counterpart of
(\ref{eq:ANStandarizedEstimator}) in which the unknown matrix $\mathbf{\Sigma
}_{n}$ is replaced by the estimator $\mathbf{\hat{\Sigma}}_{n},$ requires
additional assumptions. One possibility is to impose homoskedasticity.

\begin{theorem}
\label{thm:2}Suppose the assumptions of Theorem \ref{thm:1} hold. If
$\mathbb{E}[U_{i,n}^{2}|\mathcal{X}_{n},\mathcal{W}_{n}]=\sigma_{n}^{2},$ then
(\ref{eq:ANStudentizedEstimator}) holds with $\mathbf{\hat{\Sigma}}%
_{n}=\mathbf{\hat{\Sigma}}_{n}^{\mathtt{HO}}.$
\end{theorem}

This result shows in quite some generality that homoskedastic inference in
linear models remains valid even when $K_{n}$ is proportional to $n,$ provided
the variance estimator incorporates a degrees-of-freedom correction, as
$\mathbf{\hat{\Sigma}}_{n}^{\mathtt{HO}}$ does.

Establishing (\ref{eq:ANStudentizedEstimator}) is also possible when $K_{n}$
is assumed to be a vanishing fraction of $n,$ as is of course the case in the
usual fixed-$K_{n}$ linear regression model setup. The following theorem
establishes consistency of the conventional standard error estimator
$\mathbf{\hat{\Sigma}}_{n}^{\mathtt{EW}}$ under the assumption $\mathcal{M}%
_{n}\rightarrow_{p}0$, and also derives an asymptotic representation for
estimators of the form $\mathbf{\hat{\Sigma}}_{n}(\mathbf{\kappa}_{n})$
without imposing this assumption, which is useful to study the asymptotic
properties of other members of the HC$k$ class of standard error estimators.

\begin{theorem}
\label{thm:3}Suppose the assumptions of Theorem \ref{thm:1} hold.\newline%
(a)\ If $\mathcal{M}_{n}\rightarrow_{p}0,$ then
(\ref{eq:ANStudentizedEstimator}) holds with $\mathbf{\hat{\Sigma}}%
_{n}=\mathbf{\hat{\Sigma}}_{n}^{\mathtt{EW}}.$\newline(b) If $\Vert
\mathbf{\kappa}_{n}\Vert_{\infty}=\max_{1\leq i\leq n}\sum_{j=1}^{n}%
|\kappa_{ij,n}|=O_{p}(1),$ then%
\[
\mathbf{\hat{\Sigma}}_{n}(\mathbf{\kappa}_{n})=\frac{1}{n}\sum_{i=1}^{n}%
\sum_{j=1}^{n}\sum_{k=1}^{n}\kappa_{ik,n}M_{kj,n}^{2}\mathbf{\hat{v}}%
_{i,n}\mathbf{\hat{v}}_{i,n}^{\prime}\mathbb{E}[U_{j,n}^{2}|\mathcal{X}%
_{n},\mathcal{W}_{n}]+o_{p}(1).
\]

\end{theorem}

The conclusion of part (a) typically fails when the condition $K_{n}%
/n\rightarrow0$ is dropped. For example, when specialized to $\mathbf{\kappa
}_{n}=\mathbf{I}_{n}$ part (b) implies that in the homoskedastic case (i.e.,
when the assumptions of Theorem \ref{thm:2} are satisfied) the standard
estimator $\mathbf{\hat{\Sigma}}_{n}^{\mathtt{EW}}$ is asymptotically downward
biased in general (unless $K_{n}/n\rightarrow0$). In the following section we
make this result precise and discuss similar results for other popular
variants of the HC$k$ standard error estimators mentioned above.

On the other hand, because $\sum_{1\leq k\leq n}\kappa_{ik,n}^{\mathtt{HC}%
}M_{kj,n}^{2}=%
\I
(i=j)$ by construction, part (b) implies that $\mathbf{\hat{\Sigma}}%
_{n}^{\mathtt{HC}}$ is consistent provided $\Vert\mathbf{\kappa}%
_{n}^{\mathtt{HC}}\Vert_{\infty}=O_{p}(1).$ A simple condition for this to
occur can be stated in terms of $\mathcal{M}_{n}.$ Indeed, if $\mathcal{M}%
_{n}<1/2,$ then $\mathbf{\kappa}_{n}^{\mathtt{HC}}$ is diagonally dominant and
it follows from Theorem 1 of \cite{Varah_1975_LAA} that%
\[
\Vert\mathbf{\kappa}_{n}^{\mathtt{HC}}\Vert_{\infty}\leq\frac{1}%
{1/2-\mathcal{M}_{n}}.
\]
As a consequence, we obtain the following theorem, whose conditions can hold
even if $K_{n}/n\nrightarrow0$.

\begin{theorem}
\label{thm:4}Suppose the assumptions of Theorem \ref{thm:1} hold.\newline If
$\mathbb{P}[\mathcal{M}_{n}<1/2]\rightarrow1$ and if $1/(1/2-\mathcal{M}%
_{n})=O_{p}(1)$, then (\ref{eq:ANStudentizedEstimator}) holds with
$\mathbf{\hat{\Sigma}}_{n}=\mathbf{\hat{\Sigma}}_{n}^{\mathtt{HC}}.$
\end{theorem}

Because $\mathcal{M}_{n}\geq K_{n}/n,$ a necessary condition for Theorem
\ref{thm:4} to be applicable is that $\overline{\lim}_{n\rightarrow\infty
}K_{n}/n<1/2.$ When the design is balanced, that is, when $M_{11,n}%
=\ldots=M_{nn,n}$ (as occurs in the panel data model (\ref{eq:PanelDataModel}%
)), the condition $\overline{\lim}_{n\rightarrow\infty}K_{n}/n<1/2$ is also
sufficient, but in general it seems difficult to formulate primitive
sufficient conditions for the assumption made about $\mathcal{M}_{n}$ in
Theorem \ref{thm:4}. In practice, the fact that $\mathcal{M}_{n}$ is observed
means that the condition $\mathcal{M}_{n}<1/2$ is verifiable, and therefore
unless $\mathcal{M}_{n}$ is found to be \textquotedblleft
close\textquotedblright\ to $1/2$ there is reason to expect $\hat{\Sigma}%
_{n}^{\mathtt{HC}}$ to perform well.

\begin{description}
\item[Remark 3.] Our main results for linear models concern large-sample
approximations for the finite-sample distribution of the usual $t$-statistics.
An alternative, equally automatic approach is to employ the bootstrap and
closely related resampling procedures (see, among others,
\cite{Freedman_1981_AoS}, \cite{Mammen_1993_AoS},
\cite{Goncalvez-White_2005_JASA}, \cite{Kline-Santos_2012_JoE}). Assuming
$K_{n}/n\nrightarrow0,$ \cite{Bickel-Freedman_1983_BookCh} demonstrated an
invalidity result for the bootstrap. We conjecture that similar results can be
obtained for other resampling procedures. Furthermore, we also conjecture that
employing appropriate resampling methods on the \textquotedblleft
bias-corrected\textquotedblright\ residuals $\tilde{u}_{i,n}^{2}$ (Remark 1)
can lead to valid inference procedures. Investigating these conjectures,
however, is beyond the scope of this paper. Following the recommendation of a
reviewer, we explored the numerical performance of the standard nonparametric
bootstrap in our simulation study, where we found that indeed bootstrap
validity seems to fail in the high-dimensional settings we considered.
\end{description}

\subsection{HC$k$ Standard Errors with Many Covariates\label{section:HCk}}

The HC$k$ variance estimators are very popular in empirical work, and in our
context are of the form $\mathbf{\hat{\Sigma}}_{n}(\mathbf{\kappa}_{n})$ with
$\kappa_{ij,n}=%
\I
(i=j)\Upsilon_{i,n}M_{ii,n}^{-\xi_{i,n}}$ for some choice of $(\Upsilon
_{i,n},\xi_{i,n})$. See \cite{Long-Ervin_2000_AS} and
\cite{MacKinnon_2012_BookCh} for reviews. Theorem \ref{thm:3}(b) can be used
to formulate conditions, including $K_{n}/n\rightarrow0$, under which these
estimators are consistent in the sense that%
\[
\mathbf{\hat{\Sigma}}_{n}(\mathbf{\kappa}_{n})=\mathbf{\Sigma}_{n}%
+o_{p}(1),\qquad\mathbf{\Sigma}_{n}=\frac{1}{n}\sum_{i=1}^{n}\mathbf{\hat{v}%
}_{i,n}\mathbf{\hat{v}}_{i,n}^{\prime}\mathbb{E}[U_{i,n}^{2}|\mathcal{X}%
_{n},\mathcal{W}_{n}].
\]

More generally, Theorem \ref{thm:3}(b) shows that, if $\kappa_{ij,n}=%
\I
(i=j)\Upsilon_{i,n}M_{ii,n}^{-\xi_{i,n}}$, then%
\[
\mathbf{\hat{\Sigma}}_{n}(\mathbf{\kappa}_{n})=\mathbf{\bar{\Sigma}}%
_{n}(\mathbf{\kappa}_{n})+o_{p}(1),\qquad\mathbf{\bar{\Sigma}}_{n}%
(\mathbf{\kappa}_{n})=\frac{1}{n}\sum_{i=1}^{n}\sum_{j=1}^{n}\Upsilon
_{i,n}M_{ii,n}^{-\xi_{i,n}}M_{ij,n}^{2}\mathbf{\hat{v}}_{i,n}\mathbf{\hat{v}%
}_{i,n}^{\prime}\mathbb{E}[U_{j,n}^{2}|\mathcal{X}_{n},\mathcal{W}_{n}].
\]
We therefore obtain the following (mostly negative) results about the
properties of HC$k$ estimators when $K_{n}/n\nrightarrow0$, that is, when
potentially many covariates\ are included.

\begin{description}
\item[\textbf{HC}$0$:] $(\Upsilon_{i,n},\xi_{i,n})=(1,0).$ If $\mathbb{E}%
[U_{j,n}^{2}|\mathcal{X}_{n},\mathcal{W}_{n}]=\sigma_{n}^{2},$ then%
\[
\mathbf{\bar{\Sigma}}_{n}(\mathbf{\kappa}_{n})=\mathbf{\Sigma}_{n}%
-\frac{\sigma_{n}^{2}}{n}\sum_{i=1}^{n}(1-M_{ii,n})\mathbf{\hat{v}}%
_{i,n}\mathbf{\hat{v}}_{i,n}^{\prime}\leq\mathbf{\Sigma}_{n},
\]
with $n^{-1}\sum_{i=1}^{n}(1-M_{ii,n})\mathbf{\hat{v}}_{i,n}\mathbf{\hat{v}%
}_{i,n}^{\prime}\neq o_{p}(1)$ in general (unless $K_{n}/n\rightarrow0$).
Thus, $\mathbf{\hat{\Sigma}}_{n}(\mathbf{\kappa}_{n})=\mathbf{\hat{\Sigma}%
}_{n}^{\mathtt{EW}}$ is inconsistent in general. In particular, inference
based on $\mathbf{\hat{\Sigma}}_{n}^{\mathtt{EW}}$ is asymptotically liberal
(even) under homoskedasticity.

\item[\textbf{HC}$1$:] $(\Upsilon_{i,n},\xi_{i,n})=(n/(n-K_{n}),0).$ If
$\mathbb{E}[U_{j,n}^{2}|\mathcal{X}_{n},\mathcal{W}_{n}]=\sigma_{n}^{2}$ and
if $M_{11,n}=\ldots=M_{nn,n},$ then $\mathbf{\bar{\Sigma}}_{n}(\mathbf{\kappa
}_{n})=\mathbf{\Sigma}_{n},$ but in general this estimator is inconsistent
when $K_{n}/n\nrightarrow0$ (and so is any other scalar multiple of
$\mathbf{\hat{\Sigma}}_{n}^{\mathtt{EW}}$).

\item[\textbf{HC}$2$:] $(\Upsilon_{i,n},\xi_{i,n})=(1,1).$ If $\mathbb{E}%
[U_{j,n}^{2}|\mathcal{X}_{n},\mathcal{W}_{n}]=\sigma_{n}^{2},$ then
$\mathbf{\bar{\Sigma}}_{n}(\mathbf{\kappa}_{n})=\mathbf{\Sigma}_{n},$ but in
general this estimator is inconsistent under heteroskedasticity when
$K_{n}/n\nrightarrow0.$ For instance, if $d=1$ and if $\mathbb{E}[U_{j,n}%
^{2}|\mathcal{X}_{n},\mathcal{W}_{n}]=\hat{v}_{j,n}^{2},$ then%
\[
\mathbf{\bar{\Sigma}}_{n}(\mathbf{\kappa}_{n})-\mathbf{\Sigma}_{n}=\frac{1}%
{n}\sum_{i=1}^{n}\sum_{j=1}^{n}[\frac{M_{ij,n}^{2}}{2}(M_{ii,n}^{-1}%
+M_{jj,n}^{-1})-%
\I
(i=j)]\hat{v}_{i,n}^{2}\hat{v}_{j,n}^{2}\neq o_{p}(1)
\]
in general (unless $K_{n}/n\rightarrow0$).

\item[\textbf{HC}$3$:] $(\Upsilon_{i,n},\xi_{i,n})=(1,2).$ Inference based on
this estimator is asymptotically conservative because%
\[
\mathbf{\bar{\Sigma}}_{n}(\mathbf{\kappa}_{n})-\mathbf{\Sigma}_{n}=\frac{1}%
{n}\sum_{i=1}^{n}\sum_{j=1,j\neq i}^{n}M_{ii,n}^{-2}M_{ij,n}^{2}%
\mathbf{\hat{v}}_{i,n}\mathbf{\hat{v}}_{i,n}^{\prime}\mathbb{E}[U_{j,n}%
^{2}|\mathcal{X}_{n},\mathcal{W}_{n}]\geq0,
\]
where $n^{-1}\sum_{i=1}^{n}\sum_{j=1,j\neq i}^{n}M_{ii,n}^{-2}M_{ij,n}%
^{2}\mathbf{\hat{v}}_{i,n}\mathbf{\hat{v}}_{i,n}^{\prime}\mathbb{E}%
[U_{j,n}^{2}|\mathcal{X}_{n},\mathcal{W}_{n}]\neq o_{p}(1)$ in general (unless
$K_{n}/n\rightarrow0$).

\item[\textbf{HC}$4$:] $(\Upsilon_{i,n},\xi_{i,n})=(1,\min(4,nM_{ii,n}%
/K_{n})).$ If $M_{11,n}=\ldots=M_{nn,n}=2/3$ (as occurs when $T=3$ in the
fixed effects panel data model), then HC$4$ reduces to HC$3$, so this
estimator is also inconsistent in general.
\end{description}

Among other things these results show that (asymptotically) conservative
inference in linear models with many covariates (i.e., even when
$K/n\not \rightarrow 0$) can be conducted using standard linear methods (and
software), provided the HC$3$ standard errors are used.

In the numerical work reported in the following sections and the supplemental
appendix, we present evidence comparing all these standard error estimators.
In particular, we find that indeed standard OLS-based confidence intervals
employing HC$3$ standard errors are always quite conservative. Furthermore, we
also find that our proposed variance estimator $\mathbf{\hat{\Sigma}}%
_{n}^{\mathtt{HC}}$ delivers confidence intervals with close-to-correct
empirical coverage.

\section{Simulations\label{section:simuls}}

We conducted a simulation study to assess the finite sample properties of our
proposed inference methods as well as those of other standard inference
methods available in the literature. Based on the generic linear regression
model (\ref{eq:OLSmodel}), we consider $15$ distinct data generating processes
(DGPs) motivated by the three examples discussed above. To conserve space,
here we only discuss results from Model 1, a representative case, but the
supplemental appendix contains the full set of results and further details
(see Table 1 in the supplement for a synopsis of the DGPs used).

We discuss results for a linear model (\ref{eq:OLSmodel}) with i.i.d. data,
$n=700$, $d=1$ and $x_{i,n}\thicksim\mathsf{Normal}(0,1)$, $\mathbf{w}_{i,n}=%
\I
(\mathbf{v}_{i,n}\geq2.5)$ with $\mathbf{v}_{i,n}\thicksim\mathsf{Normal}%
(\mathbf{0},\mathbf{I}_{K_{n}})$, and $u_{i,n}\thicksim\mathsf{Normal}(0,1)$,
all independent of each other. Thus, this design considers (possibly
overlapping)\ sparse dummy variables entering $\mathbf{w}_{i,n}$; each column
assigns a value of $1$ to approximately five units out of $n=700$. We set
$\beta=1$ and $\mathbf{\gamma}_{n}=\mathbf{0}$, and considered five different
model dimensions: $\dim(\mathbf{w}_{i,n})=K_{n}\in\{1,71,141,211,281\}$. In
the supplemental appendix we also present results for more sparse dummy
variables in the context of one-way and two-way linear panel data regression
models, and for \ non-binary covariates $\mathbf{w}_{i,n}$ in both increasing
dimension linear regression settings and semiparametric partially linear
regression settings (where $\mathbf{\gamma}_{n}\not =\mathbf{0}$ and
$\mathbf{w}_{i,n}$ is constructed using power series expansions). Furthermore,
we also consider an asymmetric and a bimodal distribution for the unobservable
error terms. In all cases the numerical results are qualitatively similar to
those discussed herein. For each DGP, we investigate both homoskedastic as
well as (conditional on $x_{i,n}$ and/or $\mathbf{w}_{i,n}$) heteroskedastic
models, following closely the specifications in \cite{Stock-Watson_2008_ECMA}
and \cite{MacKinnon_2012_BookCh}. In particular, our heteroskedastic model
takes the form: $\mathbb{V}[u_{i,n}|x_{i,n},\mathbf{w}_{i,n}]=\varkappa
_{u}(1+(t(x_{i,n})+\mathbf{\iota}^{\prime}\mathbf{w}_{i,n})^{2})$ and
$\mathbb{V}[x_{i,n}|\mathbf{w}_{i,n}]=\varkappa_{v}(1+(\mathbf{\iota}^{\prime
}\mathbf{w}_{i,n})^{2})$, where the constants $\varkappa_{u}$ and
$\varkappa_{v}$ are chosen so that $\mathbb{V}[u_{i,n}]=\mathbb{V}[x_{i,n}]=1
$, and $t(a)=a%
\I
(-2\leq a\leq2)+2\operatorname*{sgn}(a)(1-%
\I
(-2\leq a\leq2))$.

We conducted $S=5,000$ simulations to study the finite sample performance of
$16$ confidence intervals: eight based on a Gaussian approximation and eight
based on a bootstrap approximation. Our paper offers theory for Gaussian-based
inference methods, but we also included bootstrap-based inference methods for
completeness (as discussed in Remark 3, the bootstrap is invalid when
$K_{n}\propto n$ in linear regression models). For each inference method, we
report both average coverage frequency and interval length of $95\%$ nominal
confidence intervals; the latter provides a summary of efficiency/power for
each inference method. To be more specific, for $\alpha=0.05$, the confidence
intervals take the form:%
\[
\mathsf{I}_{\ell}=\left[  ~\hat{\beta}_{n}-q_{\ell,1-\alpha/2}^{-1}\cdot
\sqrt{\frac{\hat{\Omega}_{n,\ell}}{n}}~,~\hat{\beta}_{n}-q_{\ell,\alpha
/2}^{-1}\cdot\sqrt{\frac{\hat{\Omega}_{n,\ell}}{n}}\right]  ,\qquad\hat
{\Omega}_{n,\ell}=\hat{\Gamma}_{n}^{-1}\hat{\Sigma}_{n,\ell}\hat{\Gamma}%
_{n}^{-1},
\]
where $q_{\ell,a}^{-1}=q_{\ell}^{-1}(a)$ and $q_{\ell}(a)$ denotes a
cumulative distribution function, and $\hat{\Sigma}_{n,\ell}$ with $\ell\in
\{$HO$0$, HO$1$, HC$0$, HC$1$, HC$2$, HC$3$, HC$4$, HC$K\}$ corresponds the
variance estimators discussed in Sections \ref{section:overview} and
\ref{section:HCk}. Gaussian-based methods set $q(a)$ equal to the standard
normal distribution for all $\ell$, while bootstrap-based methods are based on
the nonparametric bootstrap distributional approximation to the distribution
of the t-test $\mathsf{T}_{\ell}=(\hat{\beta}_{n}-\beta)/\sqrt{\hat{\Omega
}_{n,\ell}/n}$. The empirical coverage of these $16$ confidence intervals are
reported in Panel (a)\ of Table
\ref{table:simuls}%
. In addition, Panel (b)\ of Table
\ref{table:simuls}
reports the average interval length of each confidence intervals, which is
computed as $\mathsf{L}_{\ell}=(q_{\ell,1-\alpha/2}^{-1}-q_{\ell,\alpha
/2}^{-1})\cdot\sqrt{\hat{\Omega}_{n,\ell}/n}$, which offers a summary of
finite sample power/efficiency of each inference method.

The main findings from the simulation study are in line with our theoretical
results. To be precise, we find that the confidence interval estimators
constructed\ using our proposed standard errors formula $\hat{\Sigma}%
_{n}^{\mathtt{HC}}$, denoted HC$K$, offer close-to-correct empirical coverage.
The alternative heteroskedasticity consistent standard errors currently
available in the literature lead to confidence intervals that could deliver
substantial under or over coverage depending on the design and degree of
heteroskedasticity considered. We also find that inference based on HC$3$
standard errors is conservative, a general asymptotic result that is formally
established in this paper. Bootstrap-based methods seem to perform better than
their Gaussian-based counterparts, but they never outperform our proposed
Gaussian-based inference procedure nor do they provide close-to-correct
empirical coverage across all cases. Finally, our proposed confidence
intervals also exhibit very good average interval length.

\section{Empirical Illustration\label{section:empapp}}

We illustrate the different linear regression inference methods discussed in
this paper using a real data set to study the effect of ability on earnings.
In particular, we employ the dataset constructed by \cite[CHV, hereafter]%
{Carneiro-Heckman-Vytlacil_2011_AER}. [The dataset is available at
\url{https://www.aeaweb.org/articles?id=10.1257/aer.101.6.2754}%
.]. The data comes from the 1979 National Longitudinal Survey of Youth
(NLSY79), which surveys individuals born in $1957$--$1964$ and includes basic
demographic, economic and educational information for each individual. It also
includes a well-known proxy for ability (beyond schooling and work
experience): the Armed Forces Qualification Test (AFQT), which gives a measure
usually understood as a proxy for their intrinsic ability for the respondent.
This data has been used repeatedly to either control for or estimate the
effects of ability\ in empirical studies in economics and other disciplines.
See CHV for further details and references.

The sample is composed of white males of ages between $28$ and $34$ years of
old in $1991$, at most $5$ siblings, and with at least incomplete secondary
education. We split the sample into individuals with high school dropouts and
high school graduates, and individuals with some college, college graduates,
and postgraduates. For each subsample, we consider the linear regression model
(\ref{eq:OLSmodel}) with $y_{i,n}=\log(\mathtt{wages}_{i})$, where
$\mathtt{wages}_{i}$ is the log wage in 1991 of unit $i$, $x_{i,n}%
=\mathtt{afqt}_{i}$ denotes the (adjusted) standardized AFQT score for unit $i
$, and $\mathbf{w}_{i,n}$ collects several survey, geographic and dummy
variables for unit $i$. In particular, $\mathbf{w}_{i,n}$ includes the $14$
covariates described in CHV (Table 2, p. 2763), a dummy variable for wether
the education level was completed, eight cohort fixed effects, county fixed
effects, and cohort-county fixed effects. For our illustration, we further
restrict the sample to units in counties with at least $3$ survey respondents,
giving a total of $K_{n}=122$ and $n=436$ ($K_{n}/n=0.280$, $\mathcal{M}%
_{n}=0.422$) for high school educated units and $K_{n}=123$ and $n=452$
($K_{n}/n=0.272$, $\mathcal{M}_{n}=0.411$) college educated units.

The empirical findings are reported in Table
\ref{table:empapp}%
. For high school educated individuals, we find an estimated returns to
ability of $\hat{\beta}=0.060$. The statistical significance of this effect,
however, depends on the inference method employed. If homoskedastic consistent
standard errors are used, then the effect is statistical significant (p-values
are $0.010$ and $0.029$ for unadjusted and degrees-of-freedom adjusted
standard errors, respectively). If heteroskedasticity consistent standard
errors are used, the default method in most empirical studies, then the
statistical significance depends on the which inference method is used; see
Section \ref{section:HCk}. In particular, HC$0$ also gives a statistically
significant result (p-value is $0.020$), while HC$1$ and HC$2$ deliver
marginal significance (both p-values are $0.048$). On the other hand, HC3 and
HC4 give p-values of $0.092$ and $0.122$, respectively, and hence suggest that
the point estimate is not statistically distinguishable from zero. Finally,
our proposed standard error, HC$K$, gives a p-value of $0.058$, also making
$\hat{\beta}=0.060$ statistically insignificant at the conventional 5-pecent
level. In contrast, for college educated individuals, we find an effect of
$\hat{\beta}=0.091$, and all inference methods indicate that this estimated
returns to ability is statistically significant at conventional levels. In
particular, HC$3$ and our proposed standard errors HC$K$ give p-values of
$0.037$ and $0.017$, respectively.

This illustrative empirical application showcases the role of our proposed
inference method for empirical work employing linear regression with possibly
many covariates; in this application, $K_{n}$ large relative to $n$
($K_{n}/n\approx0.3$ ) is quite natural due to the presence of many county and
cohort fixed effects. Specifically, when studying the effect of ability on
earnings for high school educated individuals, the statistical significance of
the results crucially depend on the inference methods used: as predicted by
our theoretical findings, inference methods that are not robust to the
inclusion of many covariates tend to deliver statistically significant
results, while methods that are robust (HC$3$ is asymptotically conservative
and HC$K$ is asymptotically correct) do not deliver statistically significant
results, giving an example where the empirical conclusion may change depending
on whether the presence of many covariates is taken into account when
conducting inference. In contrast, the empirical findings for college educated
individuals appear to be statistically significant and robust across all
inference methods.

\section{Conclusion\label{section:conclusion}}

We established asymptotic normality of the OLS estimator of a subset of
coefficients in high-dimensional linear regression models with many nuisance
covariates, and investigated the properties of several popular
heteroskedasticity-robust standard error estimators in this high-dimensional
context. We showed that none of the usual formulas deliver consistent standard
errors when the number of covariates is not a vanishing proportion of the
sample size. We also proposed a new standard error formula that is consistent
under (conditional)\ heteroskedasticity and many covariates, which is fully
automatic and does not assume special, restrictive structure on the regressors.

Our results concern high-dimensional models where the number of covariates is
at most a non-vanishing fraction of the sample size. A quite recent related
literature concerns ultra-high-dimensional models where the number of
covariates is much larger than the sample size, but some form of
(approximate)\ sparsity is imposed in the model; see, e.g.,
\cite{Belloni-Chernozhukov-Hansen_2014_ReStud}, \cite{Farrell_2015_JoE},
\cite{Belloni-Chernozhukov-Hansen-FernandezVal_2017_ECMA}, and references
therein. In that setting, inference is conducted after covariate selection,
where the resulting number of selected covariates is at most a vanishing
fraction of the sample size (usually much smaller). An implication of the
results obtained in this paper is that the latter assumption cannot be dropped
if post covariate selection inference is based on conventional standard
errors. It would therefore be of interest to investigate whether the methods
proposed herein can be applied also for inference post covariate selection in
ultra-high-dimensional settings, which would allow for weaker forms of
sparsity because more covariates could be selected for inference.%

\onehalfspacing

\bibliographystyle{econometrica}
\bibliography{HCStdErrManyCov}

\newpage%

\begin{landscape}
\begin{table}\renewcommand{\arraystretch}{1.05}
\caption{Simulation Results (Model $1$ in Supplemental Appendix).}%
\label{table:simuls}
\subfloat[Empirical Coverage]{\resizebox{9in}{!}{
\begin{tabular}{lccccccccccccccccc}
\hline\hline
\multicolumn{1}{l}{\bfseries }&\multicolumn{8}{c}{\bfseries Gaussian Distributional Approximation}&\multicolumn{1}{c}{\bfseries }&\multicolumn{8}{c}{\bfseries Bootstrap Distributional Approximation}\tabularnewline
\cline{2-9} \cline{11-18}
\multicolumn{1}{l}{}&\multicolumn{1}{c}{HO0}&\multicolumn{1}{c}{HO1}&\multicolumn{1}{c}{HC0}&\multicolumn{1}{c}{HC1}&\multicolumn{1}{c}{HC2}&\multicolumn{1}{c}{HC3}&\multicolumn{1}{c}{HC4}&\multicolumn{1}{c}{HCK}&\multicolumn{1}{c}{}&\multicolumn{1}{c}{HO0}&\multicolumn{1}{c}{HO1}&\multicolumn{1}{c}{HC0}&\multicolumn{1}{c}{HC1}&\multicolumn{1}{c}{HC2}&\multicolumn{1}{c}{HC3}&\multicolumn{1}{c}{HC4}&\multicolumn{1}{c}{HCK}\tabularnewline
\hline
{\bfseries Homoskedastic Model}&&&&&&&&&&&&&&&&&\tabularnewline
~~$K/n=0.001$&0.949&0.950&0.948&0.948&0.948&0.948&0.948&0.948&&0.946&0.946&0.943&0.943&0.943&0.943&0.943&0.943\tabularnewline
~~$K/n=0.101$&0.939&0.956&0.939&0.952&0.952&0.962&0.980&0.951&&0.951&0.951&0.947&0.947&0.948&0.949&0.947&0.942\tabularnewline
~~$K/n=0.201$&0.916&0.947&0.919&0.947&0.946&0.968&0.989&0.945&&0.965&0.965&0.950&0.950&0.949&0.946&0.944&0.939\tabularnewline
~~$K/n=0.301$&0.900&0.950&0.904&0.954&0.951&0.977&0.983&0.949&&0.980&0.980&0.961&0.961&0.949&0.931&0.948&0.933\tabularnewline
~~$K/n=0.401$&0.881&0.954&0.884&0.955&0.952&0.989&0.972&0.949&&0.989&0.989&0.976&0.976&0.956&0.928&0.967&0.944\tabularnewline
\hline
{\bfseries Heteroskedastic Model}&&&&&&&&&&&&&&&&&\tabularnewline
~~$K/n=0.001$&0.880&0.880&0.945&0.945&0.945&0.945&0.946&0.945&&0.939&0.939&0.937&0.937&0.937&0.937&0.937&0.937\tabularnewline
~~$K/n=0.101$&0.725&0.750&0.885&0.904&0.926&0.957&0.989&0.948&&0.897&0.897&0.916&0.906&0.907&0.909&0.902&0.919\tabularnewline
~~$K/n=0.201$&0.762&0.804&0.853&0.901&0.924&0.973&0.995&0.945&&0.919&0.919&0.919&0.909&0.908&0.907&0.908&0.920\tabularnewline
~~$K/n=0.301$&0.784&0.856&0.837&0.903&0.926&0.981&0.977&0.947&&0.944&0.944&0.936&0.926&0.919&0.903&0.920&0.920\tabularnewline
~~$K/n=0.401$&0.758&0.875&0.792&0.908&0.929&0.990&0.950&0.948&&0.975&0.975&0.962&0.962&0.936&0.900&0.953&0.926\tabularnewline
\hline
\end{tabular}
}}\\
\subfloat[Interval Length]{\resizebox{9in}{!}{
\begin{tabular}{lccccccccccccccccc}
\hline\hline
\multicolumn{1}{l}{\bfseries }&\multicolumn{8}{c}{\bfseries Gaussian Distributional Approximation}&\multicolumn{1}{c}{\bfseries }&\multicolumn{8}{c}{\bfseries Bootstrap Distributional Approximation}\tabularnewline
\cline{2-9} \cline{11-18}
\multicolumn{1}{l}{}&\multicolumn{1}{c}{HO0}&\multicolumn{1}{c}{HO1}&\multicolumn{1}{c}{HC0}&\multicolumn{1}{c}{HC1}&\multicolumn{1}{c}{HC2}&\multicolumn{1}{c}{HC3}&\multicolumn{1}{c}{HC4}&\multicolumn{1}{c}{HCK}&\multicolumn{1}{c}{}&\multicolumn{1}{c}{HO0}&\multicolumn{1}{c}{HO1}&\multicolumn{1}{c}{HC0}&\multicolumn{1}{c}{HC1}&\multicolumn{1}{c}{HC2}&\multicolumn{1}{c}{HC3}&\multicolumn{1}{c}{HC4}&\multicolumn{1}{c}{HCK}\tabularnewline
\hline
{\bfseries Homoskedastic Model}&&&&&&&&&&&&&&&&&\tabularnewline
~~$K/n=0.001$&0.148&0.148&0.148&0.148&0.148&0.148&0.148&0.148&&0.148&0.148&0.149&0.149&0.149&0.149&0.149&0.149\tabularnewline
~~$K/n=0.101$&0.148&0.156&0.148&0.157&0.156&0.165&0.186&0.156&&0.161&0.161&0.158&0.158&0.158&0.158&0.158&0.157\tabularnewline
~~$K/n=0.201$&0.148&0.166&0.149&0.167&0.165&0.185&0.225&0.165&&0.180&0.180&0.170&0.170&0.169&0.167&0.166&0.164\tabularnewline
~~$K/n=0.301$&0.148&0.177&0.150&0.179&0.177&0.212&0.219&0.177&&0.210&0.210&0.189&0.189&0.182&0.172&0.180&0.174\tabularnewline
~~$K/n=0.401$&0.148&0.192&0.150&0.194&0.191&0.247&0.213&0.190&&0.260&0.260&0.223&0.223&0.200&0.174&0.212&0.189\tabularnewline
\hline
{\bfseries Heteroskedastic Model}&&&&&&&&&&&&&&&&&\tabularnewline
~~$K/n=0.001$&0.148&0.148&0.186&0.186&0.186&0.186&0.187&0.186&&0.186&0.186&0.188&0.188&0.188&0.188&0.188&0.188\tabularnewline
~~$K/n=0.101$&0.148&0.156&0.213&0.225&0.241&0.273&0.357&0.254&&0.243&0.243&0.264&0.264&0.266&0.268&0.273&0.269\tabularnewline
~~$K/n=0.201$&0.148&0.166&0.187&0.209&0.226&0.276&0.353&0.244&&0.243&0.243&0.252&0.252&0.251&0.248&0.251&0.249\tabularnewline
~~$K/n=0.301$&0.148&0.177&0.170&0.203&0.219&0.287&0.278&0.240&&0.259&0.259&0.254&0.254&0.244&0.232&0.247&0.239\tabularnewline
~~$K/n=0.401$&0.148&0.191&0.159&0.206&0.220&0.310&0.239&0.241&&0.300&0.300&0.276&0.276&0.248&0.218&0.269&0.243\tabularnewline
\hline
\end{tabular}
}}\\
\scriptsize\hspace{.5in}\raggedright Notes:
(i) DGP is Model 1 from the supplemental appendix, sample size is $n=700$, number of bootstrap replications is $B=500$, and number of simulation replications is $S=5,000$;
(ii) Columns HO$0$ and HO$1$ correspond to confidence intervals using homoskedasticity consistent standard errors without and with degrees of freedom correction, respectively,
columns HC$0$--HC$4$ correspond to confidence intervals using the heteroskedasticity consistent standard errors discussed in Sections \ref
{section:overview} and \ref{section:HCk},
and columns HC$K$ correspond to confidence intervals using our proposed standard errors estimator.
\end{table}\clearpage\end{landscape}%
%

\begin{table}\renewcommand{\arraystretch}{1.05}
\caption{Empirical Application (Returns to Ability, AFQT Score).}%
\label{table:empapp}\centering\subfloat[Secondary Education]{\resizebox{!}%
{!}{
\begin{tabular}{lcc}
\hline\hline
\multicolumn{1}{l}{}&\multicolumn{1}{c}{Outcome:}&\multicolumn{1}{c}{log(wages)}\tabularnewline
\hline
&&\tabularnewline
$\hat{\beta}$&0.060&\tabularnewline
&&\tabularnewline
&\textbf{Std.Err.}&\textbf{p-value}\tabularnewline
HO$0$&0.023&0.010\tabularnewline
HO$1$&0.028&0.029\tabularnewline
HC$0$&0.026&0.020\tabularnewline
HC$1$&0.030&0.048\tabularnewline
HC$2$&0.030&0.048\tabularnewline
HC$3$&0.036&0.092\tabularnewline
HC$4$&0.039&0.122\tabularnewline
HC$K$&0.032&0.058\tabularnewline
&&\tabularnewline
$K_n$&122&\tabularnewline
$n$&436&\tabularnewline
$K_n/n$&0.280&\tabularnewline
$\mathcal{M}_n$&0.422&\tabularnewline
\hline
\end{tabular}
}}\qquad
\subfloat[College Education]{\resizebox{!}{!}{
\begin{tabular}{lcc}
\hline\hline
\multicolumn{1}{l}{}&\multicolumn{1}{c}{Outcome:}&\multicolumn{1}{c}{log(wages)}\tabularnewline
\hline
&&\tabularnewline
$\hat{\beta}$&0.091&\tabularnewline
&&\tabularnewline
&\textbf{Std.Err.}&\textbf{p-value}\tabularnewline
HO$0$&0.032&0.005\tabularnewline
HO$1$&0.038&0.016\tabularnewline
HC$0$&0.033&0.006\tabularnewline
HC$1$&0.039&0.018\tabularnewline
HC$2$&0.038&0.016\tabularnewline
HC$3$&0.044&0.037\tabularnewline
HC$4$&0.048&0.058\tabularnewline
HC$K$&0.038&0.017\tabularnewline
&&\tabularnewline
$K_n$&123&\tabularnewline
$n$&452&\tabularnewline
$K_n/n$&0.272&\tabularnewline
$\mathcal{M}_n$&0.411&\tabularnewline
\hline
\end{tabular}
}}\\
\end{table}\clearpage

\end{document}